\documentclass[a4paper,11pt]{article}

\usepackage{graphicx}
\graphicspath{ {./Figures_Africomp/} }
\usepackage{bm}

\usepackage{authblk}
\usepackage{amscd}
\usepackage{amsfonts}
\usepackage[numbers]{natbib}
\usepackage{amsthm}
\usepackage{amssymb}
\usepackage{amsmath}
\usepackage{algpseudocode}
\usepackage{algorithm}
\usepackage{cancel}
\usepackage{chngcntr}
\usepackage{caption}
\usepackage{enumerate}
\usepackage{epsfig}
\usepackage{epic}
\usepackage{eepic}
\usepackage{enumitem}
\usepackage{etoolbox}
\usepackage{fancyhdr}
\usepackage{fullpage}
\usepackage[a4paper, total={6in, 8in}]{geometry}
\usepackage{graphicx}
\usepackage{here}
\usepackage[utf8]{inputenc}
\usepackage{lscape}
\usepackage{listings}
\usepackage{longtable}
\usepackage{natbib}
\usepackage{nccmath}
\usepackage{parskip}
\usepackage{rotating}
\usepackage{subfigure}
\usepackage{tabularx}
\usepackage{textcomp}
\usepackage{tikz}
\usepackage{titlesec}
\usepackage{threeparttable}
\usepackage{varwidth}
\usepackage{framed,multirow}

\usepackage{xcolor}

\title{Comparison of two aspects of a PDE model for biological network formation}

\author[1]{Clarissa Astuto}
\author[1]{Daniele Boffi}
\author[1]{Jan Haskovec}
\author[1]{Peter Markowich}
\author[2]{Giovanni Russo}

\affil[1]{King Abdullah University of Science and Technology (KAUST), 4700, Thuwal, Saudi Arabia}
\affil[2]{Department of Mathematics and Computer Science, University of Catania, Viale Andrea Doria 6, 95125, Catania, Italy}

\begin{document}
\maketitle

\begin{abstract}
	We compare the solutions of two systems of partial differential equations (PDE), seen as two different interpretations of the same model that describes formation of complex biological networks. Both approaches take into account the time evolution of the medium flowing through the network, and we compute the solution of an elliptic-parabolic PDE system for the conductivity vector $m$, the conductivity tensor $\mathbb{C}$ and the pressure $p$. We use finite differences schemes in a uniform Cartesian grid in the spatially two-dimensional setting to solve the two systems, where the parabolic equation is solved by a semi-implicit scheme in time. Since the conductivity vector and tensor appear also in the Poisson equation for the pressure $p$, the elliptic equation depends implicitly on time. For this reason we compute the solution of three linear systems in the case of the conductivity vector $m\in\mathbb{R}^2$, and four linear systems in the case of the symmetric conductivity  tensor $\mathbb{C}\in\mathbb{R}^{2\times 2}$, at each time step. To accelerate the simulations, we make use of the Alternating Direction Implicit (ADI) method.
	
	The role of the parameters is important for obtaining detailed solutions. We provide numerous tests with various values of the parameters involved, to see the differences in the solutions of the two systems.	
\end{abstract}	
	
\section{Introduction}
We study two elliptic-parabolic systems of partial differential equations (PDE) describing the formation of biological network structures.
Both systems are derived as gradient flows of an energy functional,
consisting of a diffusive term, activation term and metabolic cost term.
The energy functional can be seen as a continuum version of its discrete counterpart introduced in~\cite{hu2013adaptation}. Here the authors consider a total energy consumption function for a general class of biological transport networks (seen also in~\cite{katifori2010damage}), including a material cost function of the network. In the papers~\cite{hu2013adaptation,hu2013optimization}, the authors adapt the dynamics of local information for biological transport network and their relation with the optimization principle, which is known to be very common in nature.

Assuming the validity of Darcy’s law for slow flow in porous media
(see, for instance, ~\cite{darcy1856fontaines,neuman1977theoretical}), the energy functional is constrained by a Poisson equation for the fluid pressure.
We use two modes of description of the network conductivity: first, in terms
of a conductance vector $m$, and, second, in terms of a symmetric positive definite conductance tensor $\mathbb{C}$. 
Taking the $L^2$-gradient flow with respect to $m$ (see, for instance,~\cite{marko_perthame,marko_perthame_2,marko_perthame_schlo,marko_albi,portaro_arxiv}) and, resp., with respect to $\mathbb{C}$ (see~\cite{marko_pilli}),
leads to two structurally similar elliptic-parabolic PDE systems.

In the parabolic equation, a reaction term appears, representing the metabolic cost of the network, while the diffusion term describes the randomness in the material structure. The third term, called activation term, describes the tendency of the network to align with the principal direction of the material flow. The elliptic Poisson equation describes local mass conservation and is equipped with a right-hand side describing the distribution of sources and sinks of the material. This distribution is supplied as a datum and is supposed to be time independent.

The aim of this paper is to present several numerical simulations of both the vector-valued $m$-model and the tensor-valued $\mathbb{C}$-model, with the goal of comparing their solutions in various parameter settings. The initial datum is chosen such that $\mathbb{C}= m\otimes m$, i.e., initially, the principal direction (eigenvector corresponding to the largest eigenvalue) of $\mathbb{C}$ is aligned with $m$. Note that, in the spatially two-dimensional setting, the second eigenvalue of $m\otimes m$ is zero, with eigenspace orthogonal to $m$. We shall observe that the two PDE systems develop structurally and qualitatively similar solutions, however, they differ in quantitative details, for instance, the number and location of branches.

\section{Mathematical model}
We introduce the energy functional $\mathcal{E}_{\mathrm{tens}}$ for the tensor-valued model,
\begin{equation} \label{energy_T}
	\mathcal{E}_{\mathrm{tens}}[\mathbb{C}] :=
	\int_{\Omega } \frac{D^2}{2} |\nabla \mathbb{C}|^2 + c^2 \nabla p[\mathbb{C}] \cdot \mathbb{P}[\mathbb{C}] \nabla p[\mathbb{C}]  + M(|\mathbb{C}| ) \mathrm{d} x,
\end{equation}
where $\mathbb{C} = \mathbb{C}(x)\in\mathbb{R}^{2\times 2}$ is the conductivity tensor,
the diffusivity parameter $D \in \mathbb{R}$ measures the effect of random fluctuations in the network, and the activation parameter $c^2>0$ controls the strength of the network formation feedback loop.
The total permeability tensor is of the form $\mathbb{P}[\mathbb{C}] := r\mathbb{I} + \mathbb{C}$, where
the scalar function $r=r(x) \geq r_0 > 0$ describes the isotropic background permeability of the medium.
The scalar pressure $p=p[\mathbb{C}]$ of the fluid transported within the network is the unique solution (up to an additive constant) of the Poisson equation
\begin{equation}
	- \nabla\cdot \left( \mathbb{P}[\mathbb{C}] \nabla p \right) = S, \label{eq:Poisson}
\end{equation}
subject to homogeneous Neumann boundary condition on $\partial\Omega$.
The source/sink distribution $S=S(x)$ in the mass conservation Eq. \eqref{eq:Poisson} is to be supplemented as an input datum and is assumed to be independent of time.
The metabolic cost function $M:\mathbb{R}^+ \to \mathbb{R}^+$ describes the dependence of the metabolic expenditure of maintaining the network on its transportation capacity, see~\cite{hu2013adaptation}.
The expression $|\mathbb{C}|$ denotes the Frobenius norm $|\mathbb{C}| := \sqrt{\sum_{i=1}^2 \sum_{j=1}^2 {C}_{ij}^2}$ and $|\nabla \mathbb{C}| := \sqrt{\sum_{i=1}^2 \sum_{j=1}^2 \sum_{k=1}^2 (\partial{C}_{ij}/\partial x_k)^2}$.

Taking the $L^2$-gradient flow of the energy \eqref{energy_T} constrained by \eqref{eq:Poisson}, we obtain the parabolic-elliptic system
\begin{align}
	\label{eq_darcy_p_tens}
	-\nabla\cdot \left( \left(r\mathbb{I} + \mathbb{C} \right) \nabla p \right) &= S \\
	\label{eq_reaction_diff_tens}
	\frac{\partial\mathbb{C}}{\partial t} - D^2\Delta \mathbb{C} -c^2\nabla p \otimes \nabla p + \frac{M'\left(|\mathbb{C}|\right)}{\left| \mathbb{C} \right|} \mathbb{C} & = 0,
\end{align}
see~\cite{marko_pilli} for details of the derivation. In our paper, we shall make the generic choice for the metabolic cost function (see, e.g.,~\cite{marko_albi,marko_perthame,marko_perthame_2,marko_perthame_schlo}),
\begin{equation} \label{M}
	M(s) := \frac{\alpha}{\gamma} s^\gamma\qquad\mbox{for } s\geq 0,
\end{equation}
where $\alpha>0$ is the metabolic constant and $\gamma>0$ the metabolic exponent.
For instance, to model leaf venation in plants, one chooses $1/2<\gamma<1$, see~\cite{hu2013adaptation}. Then Eq.~\eqref{eq_reaction_diff_tens} becomes
\begin{equation}
	\label{eq_C_metab}
	\frac{\partial\mathbb{C}}{\partial t} - D^2\Delta \mathbb{C} -c^2\nabla p \otimes \nabla p + \alpha |\mathbb{C}|^{\gamma-2} \mathbb{C} = 0  
\end{equation} 

To derive the vector-valued model, we make the ansatz $\mathbb{C}:=m\otimes m$. Inserting this into \eqref{energy_T} with $D=0$ and noting that $|\mathbb{C}|=|m|^2$, we obtain
\begin{equation*}
	\mathcal{E} = \int_{\Omega } c^2 \nabla p[m] \cdot \mathbb{P}[m] \nabla p[m]  + M\left(|m|^2\right) \mathrm{d} x
\end{equation*}
with $\mathbb{P}[m] = r\mathbb{I} + m\otimes m$ and, with a slight abuse of notation, we now denote by $p=p[m]$ the solution of the Poisson equation \eqref{eq:Poisson} with $\mathbb{P}[\mathbb{C}]$ replaced by $\mathbb{P}[m\otimes m]$.
Re-introducing the Dirichlet integral $D^2 \int_\Omega |\nabla m|^2 \mathrm{d} x$, we arrive at the energy functional
\begin{equation} \label{energy_m}
	\mathcal{E}_{\mathrm{vect}}[m] :=
	\int_{\Omega } D^2 |\nabla m|^2 + c^2 \nabla p[m] \cdot \mathbb{P}[m\otimes m] \nabla p[m]  + M(|m|^2) \mathrm{d} x.
\end{equation}
Note that this functional is different from the one obtained by replacing $C= m\otimes m$ in Eq.~\eqref{energy_T}, the two functionals mainly differencing in the contribution of the diffusion term. As we shall see, we are mainly interested in studying the behavior of the system for very small diffusion coefficients, so that we expect such difference not being so crucial.

Taking the $L^2$-gradient flow with respect to the vector variable $m$ and using the explicit form \eqref{M} for the metabolic function, we obtain the system
\begin{align}
	\label{sy_m_p}
	-\nabla\cdot \left( \left(r\mathbb{I} + m \otimes m \right) \nabla p \right) &= S \\
	\label{sy_m_diff}
	\frac{\partial {m}}{\partial t} - D^2\Delta m - c^2(\nabla p \otimes \nabla p) m + {\alpha}|m|^{2(\gamma-1)} m & = 0.
\end{align}

From now on, we shall denote system~\eqref{eq_darcy_p_tens}--\eqref{eq_C_metab} as the $\mathbb{C}-$system, while we call \eqref{sy_m_p}--\eqref{sy_m_diff} the $m-$system.
Let us now shortly discuss the differences between the two systems.
First, the activation term $c^2\nabla p\otimes\nabla p$ in \eqref{eq_reaction_diff_tens} is quadratic and depends on $\mathbb{C}$ only through the solution $p=p[\mathbb{C}]$ of the Poisson equation \eqref{eq_darcy_p_tens}. In contrast, the activation term $c^2 (\nabla p\otimes\nabla p)m$ in \eqref{sy_m_diff} is cubic.
Moreover, the metabolic term $\alpha |\mathbb{C}|^{\gamma-2}\mathbb{C}$ in \eqref{eq_reaction_diff_tens} becomes singular at $\mathbb{C}=0$ if (and only if) $\gamma<1$. This obviously causes difficulties for both the analytical treatment of the system (see~\cite{marko_pilli} for details) and its numerical resolution. We shall discuss in Section \ref{section_C_time_discr} how the numerical difficulties can be overcome. In contrast, the metabolic term ${\alpha}|m|^{2(\gamma-1)} m$ in \eqref{sy_m_diff} only becomes singular at $m=0$ if $\gamma<1/2$. Taking into account that for typical applications in biology only the parameter range $\gamma\geq 1/2$ is relevant (see~\cite{hu2013adaptation}), the metabolic term in \eqref{sy_m_diff}) does not require any special treatment for the $m-$system.

We define the equations on a domain $\Omega \subset \mathbb{R}^2$, and we define the boundary conditions on $\partial \Omega$ for the pressure and the conductivity. We choose homogeneous Dirichlet boundary conditions for $m$ and $\mathbb{C}$, and homogeneous Neumann conditions for $p$. 
\begin{equation}
	\label{eq_bc}
	m(t,\vec{x}) = 0, \quad \mathbb{C}(t,\vec{x}) = 0, \quad \mathbb{P}[\mathbb{C}]\nabla p(t,\vec{x}) \cdot \nu = 0, \quad \vec{x}\in \partial \Omega, \, t\geq 0
\end{equation}
where $\nu$ is the outgoing normal vector to $\partial \Omega$.
In all our numerical test, the numerical support of $\mathbb{C}$ or $m$, after a long time, is well within $\Omega$, therefore the solution should not be particularly sensitive to the boundary conditions on the two variables. In any case, different boundary conditions for $m$ and $\mathbb{C}$ are currently under investigation.

To close the system, we prescribe an initial condition for the conductivity vector and tensor
\begin{equation}
	\label{eq_ic}
	m(t=0,\vec{x}) = m^0(\vec{x}), \quad \mathbb{C}(t=0,\vec{x}) = \mathbb{C}^0(\vec{x}), \quad \text{ in } \Omega.
\end{equation}
A direct consequence of the boundary condition defined for the pressure, that is also a necessary condition for the solvability of the Poisson equation, is that the source function has to be vanishing mean, i.e., $\displaystyle \int_\Omega S(\vec{x}) d\Omega = 0$.

\section{Numerical schemes}
\label{section_methods}
In this section we define the numerical schemes used to discretize in space and time the Eqs.~(\ref{sy_m_p} - \ref{sy_m_diff}) and (\ref{eq_darcy_p_tens}-\ref{eq_reaction_diff_tens}). We adopt semi-implicit second order schemes.

\subsection{Space discretization}
\label{section_space_discrete}
In space we consider the square domain $\Omega = [0,1]\times [0,1]$, that we discretize by a uniform Cartesian mesh with spatial step $h := \Delta x = \Delta y$. We call $\Omega_h$ the discrete computational domain. The two variables conductivity and the pressure are ${ \mathbb{C}_{ij}\approx \mathbb{C}(x_i,y_j)}, m_{ij} \approx m(x_i,y_j)$ and ${ p_{ij}\approx p(x_i,y_j)}$, defined at the center of the cell $(i,j)$, therefore the set of grid points is $x_i=(i-1/2)h, \, y_j=(j-1/2)h, \,(i,j)\in\{1,\dots,N\}^2$, $h N = 1$. 

In order to obtain second order accuracy in space, we use central differences for the computation of the {space derivatives}~\cite{Wesseling2023600}. 
Discretizing in space Eq.~\eqref{sy_m_diff}, we have:
\begin{eqnarray}
	\label{eq_m1_sp}	\displaystyle \frac{\partial m^{(1)}}{\partial t}&=& D^2\mathcal{L}\,m^{(1)} + c^2\left(\mathcal{D}_x p\right)^2\, m^{(1)} + c^2\mathcal{D}_x p \mathcal{D}_y p\, m^{(2)}  - \alpha |m|^{2(\gamma - 1)}m^{(1)} \\ 
	\label{eq_m2_sp}	\displaystyle \frac{\partial m^{(2)}}{\partial t}&=& D^2\mathcal{L}\,m^{(2)} + c^2\left(\mathcal{D}_y p\right)^2\, m^{(2)} + c^2\mathcal{D}_x p \mathcal{D}_y p\, m^{(1)} -  \alpha |m|^{2(\gamma - 1)}m^{(2)}
\end{eqnarray}
while the semi-discrete version of Eq.~\eqref{eq_reaction_diff_tens} is
\begin{eqnarray}
	\label{eq_C11}	\displaystyle \frac{\partial {C}^{(1,1)}}{\partial t}&=& D^2\mathcal{L}\,{C}^{(1,1)} + c^2\left(\mathcal{D}_x p\right)^2 - \gamma |\mathbb{C}|^{\gamma - 2}{C}^{(1,1)} \\ 
	\label{eq_C12}
	\displaystyle \frac{\partial {C}^{(1,2)}}{\partial t}&=& D^2\mathcal{L}\,{C}^{(1,2)} + c^2\mathcal{D}_x p \mathcal{D}_y p - \gamma |\mathbb{C}|^{\gamma - 2}{C}^{(1,2)} \\ \label{eq_C22}		\displaystyle \frac{\partial {C}^{(2,2)}}{\partial t}&=& D^2\mathcal{L}\,{C}^{(2,2)} + c^2\left(\mathcal{D}_y p\right)^2 - \gamma |\mathbb{C}|^{\gamma - 2}{C}^{(2,2)}
\end{eqnarray}
where $\mathcal{L}$ is the discrete Laplacian operator and $\mathcal{D}_x$ and $\mathcal{D}_y$ are, respectively, the discrete $x$ and $y$ first derivative operators, both using central difference approximation. The norm $|\cdot|$ we use is the Frobenius one. In this case we have three equations instead of four, because $\mathbb{C}$ is symmetric.

In order to have a fully implicit scheme, we define a compact form of the semi-discrete Eqs.~(\ref{eq_m1_sp}-\ref{eq_m2_sp}) and Eqs.~(\ref{eq_C11}-\ref{eq_C22}), as follows:
\begin{eqnarray}
	\label{eq_m_compact}
	\displaystyle \frac{\partial {m}_{\rm comp}}{\partial t}&=& D^2\mathcal{L}\,{m}_{\rm comp} + c^2 \mathcal{P}\,{m}_{\rm comp} - \alpha \mathcal{Q}^m(m){m}_{\rm comp} \\
	\label{eq_C_compact}	
	\displaystyle \frac{\partial {\mathbb{C}}_{\rm comp}}{\partial t}&=& D^2\mathcal{L}\,{\mathbb{C}}_{\rm comp} + c^2 \mathcal{P} - \alpha \mathcal{Q}^c(\mathbb{C}) {\mathbb{C}}_{\rm comp}
\end{eqnarray}
where ${m}_{\rm comp} = [m^{(1)},m^{(2)}]^T$, ${\mathbb{C}}_{\rm comp} = [C^{(1,1)},C^{(1,2)},C^{(2,2)}]^T$ and  $\mathcal{P} = \mathcal{D}_\beta p\,\mathcal{D}_\eta p $, with the suitable choice of $\beta,\eta \in \{x,y\}$. For the metabolic terms, we have
\begin{eqnarray}
	\label{eq_m_Qcal}
	\mathcal{Q}^m(m) &=&  |m|^{2(\gamma - 1)} \\
	\label{eq_Qcal}
	\mathcal{Q}^c(\mathbb{C}) &=& |\mathbb{C} |^{\gamma - 2}.
\end{eqnarray}


Here, we discretize the Poisson equation for the pressure. To obtain a conservative scheme, we consider the following discretization: first, we extend the formula in Eq.~\eqref{eq_darcy_p_tens}, that becomes
\begin{align}
	\partial_x\left(\left(r + C^{(1,1)} \right) \partial_x p \right) + \partial_x\left( C^{(1,2)} \partial_y p \right) + \partial_y\left( C^{(1,2)} \partial_x p \right)  
	+ \partial_y\left(\left(r + C^{(2,2)} \right) \partial_y p \right) = -S
\end{align}
where we use the symmetry $C^{(1,2)} = C^{(2,1)}$. 

Now, we discretize the components of the formula, one by one, since we use different discretizations. For simplicity we pose $\mathcal{C}^{1,1} = r + C^{(1,1)}$:
\begin{align} \nonumber
	\partial_x\left(\mathcal{C}^{1,1} \partial_x\, p \right)_{i,j} \approx &
	\frac{1}{\Delta x} \left( \mathcal{C}^{1,1}_{i+1/2,j} \partial_x \,p_{i+1/2,j} - \mathcal{C}^{1,1}_{i-1/2,j}\partial_x \,p_{i-1/2,j} \right)\\ \nonumber
	= & \frac{1}{2\Delta x^2} \left( \left( \mathcal{C}^{1,1}_{i+1,j} + \mathcal{C}^{1,1}_{i,j}\right)p_{i+1,j} + \left( \mathcal{C}^{1,1}_{i-1,j} + \mathcal{C}^{1,1}_{i,j}\right)p_{i-1,j} \right)\\
	& - \frac{1}{2\Delta x^2} \left( \mathcal{C}^{1,1}_{i+1,j} + \mathcal{C}^{1,1}_{i-1,j} + 2\mathcal{C}^{1,1}_{i,j}\right)p_{i,j}    
\end{align}
where in the last line we consider the following approximations:
\[ \partial_x \,p_{i+1/2,j} \approx \frac{p_{i+1,j} - p_{i,j}}{\Delta x}, \quad \mathcal{C}^{1,1}_{i+1/2,j} \approx \frac{\mathcal{C}^{1,1}_{i+1,j} + \mathcal{C}^{1,1}_{i,j}}{2}.
\]
We omit the term with both $y$-derivatives because it is analogue to the one with $x$-derivatives. Now we discretize the term with mix derivatives.
\begin{align}
	\partial_x\left(C^{(1,2)} \partial_y\, p \right)_{i,j} \approx &
	\frac{1}{\Delta x} \left( C^{(1,2)}_{i+1/2,j} \partial_y p_{i+1/2,j} - C^{(1,2)}_{i-1/2,j} \partial_y p_{i-1/2,j} \right) \\ \nonumber
	= & \frac{1}{8 \Delta x^2}\left( C^{(1,2)}_{i+1,j} + C^{(1,2)}_{i,j}\right)(p_{i+1,j+1} - p_{i+1,j-1})  
	\\ &-\frac{1}{8 \Delta x^2} \left( C^{(1,2)}_{i-1,j} + C^{(1,2)}_{i,j}\right)(p_{i-1,j+1} - p_{i-1,j-1})   \\
	& + \frac{1}{8 \Delta x^2}\left( C^{(1,2)}_{i+1,j} - C^{(1,2)}_{i-1,j} \right) (p_{i,j+1} - p_{i,j-1}) 
\end{align}
again, we omit the term with $y,x$-derivatives because it is analogue to the one with $x,y$-derivatives.

\subsection{Time discretization: \textit{symmetric}-ADI method}
At this stage, we describe the time discretization that we apply to the model. Since  systems~(\ref{sy_m_p}-\ref{eq_C_metab}) and (\ref{eq_darcy_p_tens}-\ref{eq_reaction_diff_tens}) are 
stiff in all their components, the choice of the time discretization is crucial for the efficiency. 

We also need a high performing scheme in time, since at each time step we compute the solution of seven linear systems (two for the conductivity vector $m$, three for the conductivity tensor $\mathbb{C}$ and two Poisson equations for the pressure $p$). 

As we shall see, for some values of the parameters, the well-posedness of the problem becomes weaker, which reflects the bad conditioning of the numerical problem. For such a reason, we adopt a symmetric scheme, which better preserves possible symmetries of the solution. In particular, we adopt a \textit{symmetric}-ADI scheme for both the conductivity variables, which guarantees efficiency, second order accuracy and spatial symmetry.

Anyway, the scheme is not strictly second order accurate in time for two reasons. First, the pressure is computed at time $n$ rather than at
an intermediate time $n+1/2$. Second, the metabolic term is treated partially explicitly and partially implicitly, thus destroying second order accuracy.
Improvements of the order of accuracy in time are currently under investigation.

\subsubsection{Time discretization for the conductivity vector}
\label{section_C_time_discr}
Here we focus on the time discretization for the Eqs.~(\ref{sy_m_p}-\ref{sy_m_diff}).  


Given $m^n \approx m(t^n)$, we compute $p^{n}$ by solving the Poisson equation 
\begin{equation}
	\label{eq_poisson_discrete}
	-\mathcal{L}\left(m^{n}\otimes m^n\right)\, p^{n} = S,
\end{equation}
where $\mathcal{L}\left(m^{n}\otimes m^n\right) \in \mathbb{R}^{N^2 \times N^2}$ is the discrete elliptic operator, in both directions ($x$ and $y$), with variable coefficients and corresponding to zero Neumann conditions.

The \textit{symmetric}-ADI method to solve the Eq.~\eqref{eq_m_compact} works as follows. We start with the $y$-direction implicit and $x$-direction explicit, and then we consider the opposite order in the second step of the ADI scheme 
\begin{eqnarray*}
	(y-{\rm impl},x-{\rm expl}) \qquad	\tilde{m} &=& m^n  + \frac{\Delta t}{2}\mathcal{L}_y\tilde{m} + \frac{\Delta t}{2}\mathcal{L}_x m^n + \Delta t\,c^2 \mathcal{P}_{xy}^{n}(m^{n}) 
\end{eqnarray*} 
and we solve for $\tilde{m}$,
\begin{eqnarray*}
	(x-{\rm impl},y-{\rm expl}) \qquad  m^{n+1}_y&=& \tilde{m} + \frac{\Delta t}{2}\mathcal{L}_x{ m^{n+1}_y} + \frac{\Delta t}{2}\mathcal{L}_y\tilde{m}   - \Delta t\,\alpha \mathcal{Q}^m\left({ m^{n}}\right)m^{n+1}_y \\ 
	&& + \Delta t\,c^2 \mathcal{P}_y^{n}\, m_y^{n+1}
\end{eqnarray*} 
and we solve for $m^{n+1}_y$.

The second time we apply the ADI method, we first consider the $x$-direction implicit and $y$ explicit, and then we exchange the order. Thus we have
\begin{eqnarray*}
	(x-{\rm impl},y-{\rm expl}) \qquad	\hat{m} &=& m^n + \frac{\Delta t}{2}\mathcal{L}_x\hat{m} + \frac{\Delta t}{2}\mathcal{L}_y{m}^n  + \Delta t\,c^2 \mathcal{P}_{xy}^{n}(m^{n}) 
\end{eqnarray*} 
here we solve for $\hat{m}$,
\begin{eqnarray*}
	(y-{\rm impl},x-{\rm expl}) \qquad  	m^{n+1}_x &=& \hat{m} + \frac{\Delta t}{2}\mathcal{L}_y{ m^{n+1}_x} + \frac{\Delta t}{2}\mathcal{L}_x\hat{m} - \Delta t\,\alpha \mathcal{Q}^m\left( { m^{n}}\right)m_x^{n+1} \\ 
	&& + \Delta t\,c^2 \mathcal{P}_x^{n}\, m_x^{n+1}
\end{eqnarray*} 
and now we solve for $m^{n+1}_x$, where $\mathcal{Q}^m\left(m^{n}\right) = |m^n|^{2(\gamma - 1)}$ and $ \displaystyle \mathcal{L}_\beta $, with $  \beta=x,y$, { are the discrete operators for the second derivatives} in $x$ and $y$ directions, respectively, with $  \displaystyle \mathcal{L}_\beta \in  \mathbb{R}^{N\times N}$. For the pressure term we have $\mathcal{P}_x^n = (\mathcal{D}_x p^n)^2$ for the first component of the vector $m^{(1)}$, and $\mathcal{P}_y^n = (\mathcal{D}_y p^n)^2$ for the second component $m^{(2)}$. While the force term is  $\mathcal{P}^n_{xy}(m^n)  = [\mathcal{D}_x p^n \mathcal{D}_y p^n\, m^{(2)n}, \mathcal{D}_x p^n \mathcal{D}_y p^n\, m^{(1)n}]^T$. 

At the end, we calculate the average of the two solutions $m^{n+1}_y$ and $m^{n+1}_x$ to obtain the conductivity vector at time $t^{n+1}$
\begin{eqnarray}
	\label{eq_m_ADI_final}
	m^{n+1} &=& \frac{1}{2}{ m^{n+1}_x} + \frac{1}{2}{ m^{n+1}_y}.
\end{eqnarray}

\subsubsection{Time discretization for the conductivity tensor}
As we said in the description of the two systems, for the conductivity tensor $\mathbb{C}$, the reaction term is very stiff because of the exponent $\gamma$ that belongs to the interval $(0.5,1)$. After a finite time, we are basically dividing by zero at each time step. For this reason we introduce a small regularizing parameter $\varepsilon$ in the equation, as follows
\begin{equation}
	\label{eq_Qcal2}
	\mathcal{Q}^c(\mathbb{C}) =  |\mathbb{C} + \varepsilon |^{\gamma - 2}\mathbb{C}
\end{equation}
and we study the behaviour of the system as $\varepsilon$ becomes smaller and smaller.

Given the pressure at time $t^n$ from Eq.~\eqref{eq_poisson_discrete}, we apply the \textit{symmetric}-ADI method to solve the Eq.~\eqref{eq_C_compact}. As explained before, we first choose the $y$-direction implicit and then we exchange the two directions. The scheme reads
\begin{eqnarray*}
	(y-{\rm impl},x-{\rm expl}) \qquad	\tilde{\mathbb{C}}_1 &=& \mathbb{C}^n + \frac{\Delta t}{2}\mathcal{L}_y\tilde{\mathbb{C}}_1  + \frac{\Delta t}{2}\mathcal{L}_x\mathbb{C}^n  + \Delta t\,\mathcal{P}^{n}
\end{eqnarray*} 
and we solve for $\tilde{\mathbb{C}}_1$,
\begin{eqnarray*}
	(x-{\rm impl},y-{\rm expl}) \qquad  	\mathbb{C}^{n+1}_y&=& \tilde{\mathbb{C}}_1 + \frac{\Delta t}{2}\mathcal{L}_x{ \mathbb{C}^{n+1}_y} + \frac{\Delta t}{2}\mathcal{L}_y\tilde{\mathbb{C}}_1 - \Delta t\,\alpha \mathcal{Q}^c\left({ \mathbb{C}^{n}}\right)\mathbb{C}^{n+1}_y 
\end{eqnarray*} 
and we solve for $\mathbb{C}^{n+1}_y$. 

Now, as before, we apply the ADI method for the second time for $\mathbb{C}$ starting with the $x$-direction implicit, and we have
\begin{eqnarray*}
	(x-{\rm impl},y-{\rm expl}) \qquad	\tilde{\mathbb{C}}_2 &=& \mathbb{C}^n + \frac{\Delta t}{2}\mathcal{L}_x\tilde{\mathbb{C}}_2 + \frac{\Delta t}{2}\mathcal{L}_y{\mathbb{C}}^n  + \Delta t\,\mathcal{P}^{n} 
\end{eqnarray*} 
here we solve for $\tilde{\mathbb{C}}_2$,
\begin{eqnarray*}
	(y-{\rm impl},x-{\rm expl}) \qquad  	\mathbb{C}^{n+1}_x &=& \tilde{\mathbb{C}}_2 + \frac{\Delta t}{2}\mathcal{L}_y{ \mathbb{C}^{n+1}_x} + \frac{\Delta t}{2}\mathcal{L}_x\tilde{\mathbb{C}}_2 - \Delta t\,\alpha \mathcal{Q}^c\left( { \mathbb{C}^{n}}\right)\mathbb{C}_x^{n+1}
\end{eqnarray*} 
and here we solve for $\mathbb{C}^{n+1}_x$. 

Finally we calculate $\mathbb{C}^{n+1}$ from the two solutions $\mathbb{C}^{n+1}_y$ and $\mathbb{C}^{n+1}_x$
\begin{eqnarray*}
	\mathbb{C}^{n+1} &=& \frac{1}{2}{ \mathbb{C}^{n+1}_x} + \frac{1}{2}{ \mathbb{C}^{n+1}_y}.
\end{eqnarray*}
The quantities above have the following expressions: $\mathcal{Q}^c\left(\mathbb{C}^{n}\right) = |\mathbb{C}^n + \varepsilon|^{\gamma - 2}$, and for the pressure, $\mathcal{P}^n = \mathcal{D}_\beta p^n\,\mathcal{D}_\eta p^n $, with the suitable choice of $\beta,\eta \in \{x,y\}$ for the four components of the conductivity tensor.

\section{Numerical results}
In this section we perform several simulations with the aim of studying the effect of the various parameters. In particular, we check the agreement of the two models for the $m-$system in Eqs.~(\ref{sy_m_p}-\ref{sy_m_diff}) and for the $\mathbb{C}-$system in Eqs.~(\ref{eq_darcy_p_tens}-\ref{eq_C_metab}). 

\subsection{Accuracy tests and qualitative agreements }
\label{section_qualit_test}
In Table~\ref{tab_qualit} we define the tests we want to show in this paper, varying the parameters of the systems. This choice of parameters, a typical time scale is of the order of unit, while after time 15, the solution reached the steady state.

First, we check the accuracy of the schemes adopted. In Table~\ref{table_accuracy_m} and Table~\ref{table_accuracy_C} we see the error for the conductivity variables, calculated with Richardson extrapolation (see, e.g.,~\cite{richardson1911ix}). We show the error for the module of the vector and of the tensor, and the parameters chosen are defined in \textsc{TestA}, \textsc{TestB} and \textsc{TestC}.

\begin{table}[H]
	\centering
	\begin{tabular}{|c|c|c|c|c|c|c|c|c|}  \hline
		& & $\alpha$ & $c$ & $D$ & $\varepsilon$ & $\gamma$ & $r$ & $t_{\rm fin}$ \\ \hline 
		Accuracy $m$  & \textsc{TestA}: & 0.5 & 1 & 0.01 & - & 0.75 & 0.1 & 1 \\ \hline
		Accuracy $\mathbb{C}$  &  \textsc{TestB}: &  1 & 1 & 0.01 & 0.1 & 1.75 & 0.1 & 1 \\ \hline
		Accuracy $m$ & \textsc{TestC}: &  0.5 & 5 & 0.01 & - & 0.75 & 0.01 & 1 \\ \hline \hline
		$D = 0.05$ & \textsc{TestG}: &  0.75 & 5 & 0.05 & $10^{-3}$ & 0.75 & 0.005 & 15 \\ \hline
		$D = 0.01$ & \textsc{TestD}: &  0.75 & 5 & 0.01 & $10^{-3}$ & 0.75 & 0.005 & 15 \\ \hline
		$D = 0.001$ & \textsc{TestE}: & 0.75 & 5 & 0.001 & $10^{-3}$  & 0.75 & 0.005 & 15 \\ \hline \hline
		$\gamma = 1$ & \textsc{TestH}: &  0.75 & 5 & 0.01 & $10^{-3}$ & 1 & 0.005 & 15 \\ \hline
		$\gamma = 0.75$ & \textsc{TestD}: & 0.75 & 5 & 0.01 & $10^{-3}$ & 0.75 & 0.005 & 15 \\ \hline
		$\gamma = 0.5$ & \textsc{TestF}: & 0.75 & 5 & 0.01 & $10^{-3}$ & 0.5 & 0.005 &  15 \\ \hline \hline
		$\varepsilon = 10^{-2}$ & \textsc{TestI}: &  0.75 & 5 & 0.01 & $10^{-2}$ & 0.75 & 0.005 & 15 \\ \hline
		$\varepsilon = 10^{-3}$ & \textsc{TestD}: &  0.75 & 5 & 0.01 & $10^{-3}$ & 0.75 & 0.005 & 15. \\ \hline
		$\varepsilon = 10^{-4}$ & \textsc{TestL}: &  0.75 & 5 & 0.01 & $10^{-4}$ & 0.75 & 0.005 & 15 \\ \hline
	\end{tabular}
	\caption{\textit{In this table we define all the tests that we show in Section\ref{section_qualit_test}. The first three rows show the parameters for the accuracy tests for $m$ and $\mathbb{C}$, and the results are summarized in Table~\ref{table_accuracy_m}-\ref{table_accuracy_C}. The second three rows define the parameters that we use in Fig.~\ref{figure_diff}, where we compare the results changing the diffusivity. The third three rows are the tests showed in Fig.~\ref{fig_gamma}, varying $gamma$ and the results of the last three rows are in Fig.~\ref{fig_epsilon}, where we change the stabilization parameter $\varepsilon$. For the accuracy tests, the number of points of the discretization is specified in Table~\ref{table_accuracy_m}-\ref{table_accuracy_C}, while, for all the other tests, the number of points is fixed and it is $N = 600$.}}
	\label{tab_qualit}
\end{table}

Here we define the initial conditions $m^0_{\rm comp}(\vec{x})$ and $\mathbb{C}^0_{\rm comp}(\vec{x})$, and the source function $S(\vec{x})$
\begin{eqnarray} \label{eq_ic_m}
	&m^0_{\rm comp}(\vec{x}) = [1,1]^T,\qquad \mathbb{C}^0_{\rm comp}(\vec{x
	}) = [1,0,1]^T,  \qquad S(\vec{x}) = E - \bar{E}\\
	& E = \exp(-\sigma(\vec{x}-\vec{x}_0)^2), \, \sigma = 1000, \, \vec{x}_0 = (0.1,0.1)
\end{eqnarray}
where $\mathbb{I}$ is the identity matrix and $\bar{E} = {\rm mean}(E)$.

\begin{table}[H]
	\centering
	\begin{tabular}{||c||c||c||} \hline \hline
		N & error & order \\ \hline \hline
		20 & - & - \\ \hline \hline
		40 & 0.036030 & - \\ \hline \hline
		80 & 0.0492860 & -0.4520 \\ \hline \hline
		160 & 0.01454106 & 1.7610  \\ \hline \hline
		320 & 0.00690830 & 1.0737 \\ \hline \hline
		640 & 0.001529779 &  2.1750 \\ \hline \hline
	\end{tabular}
	\begin{tabular}{||c||c||c||} \hline \hline
		N & error & order \\ \hline \hline
		20 & - & - \\ \hline \hline
		40 & 0.036012 & - \\ \hline \hline
		80 & 0.0493010 & -0.4531 \\ \hline \hline
		160 & 0.01456192 & 1.7594  \\ \hline \hline
		320 & 0.00691103 & 1.0752 \\ \hline \hline
		640 & 0.001528055 &  2.1772 \\ \hline \hline
	\end{tabular}    
	\caption{\textit{Accuracy test of the $m$-system~(\ref{sy_m_p}-\ref{sy_m_diff}): we show the $L^2-$norm of the relative error for $|m|$, with the parameters defined in} \textsc{TestA} \textit{(left) and} \textsc{TestC} \textit{(right)}.}
	\label{table_accuracy_m}
\end{table}

\begin{table}[H]
	\centering
	\begin{tabular}{||c||c||c||}
		\hline \hline
		$N$ & error$_2$ & order  \\ \hline \hline
		25 &  - & - \\ \hline \hline
		50 & $9.066 \times 10^{-2}$  & - \\ \hline \hline
		100 & $4.625 \times 10^{-2}$ & 0.97 \\ \hline \hline
		200 & $1.571 \times 10^{-2}$ & 1.56  \\ \hline \hline
		400 & $4.149 \times 10^{-3}$  & 1.92 \\ \hline \hline
		800 & $7.347 \times 10^{-4}$  & 2.50 \\ \hline \hline
	\end{tabular}
	\caption{\textit{Accuracy test of the $\mathbb{C}$-system~(\ref{eq_darcy_p_tens}-\ref{eq_reaction_diff_tens}): we show the $L^2-$norm of the relative error for $|\mathbb{C}|$, with the parameters defined in} \textsc{TestB}.}
	\label{table_accuracy_C}
\end{table}

\begin{figure}[H]
	\centering
	\begin{minipage}{.32\textwidth}    \includegraphics[width=1.\textwidth]{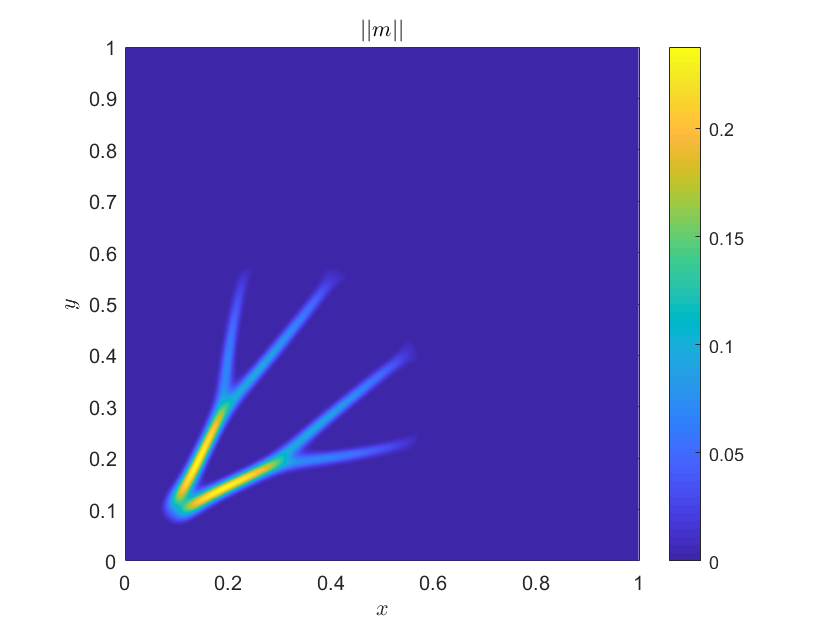}
	\end{minipage}
	\begin{minipage}{.32\textwidth}    \includegraphics[width=1\textwidth]{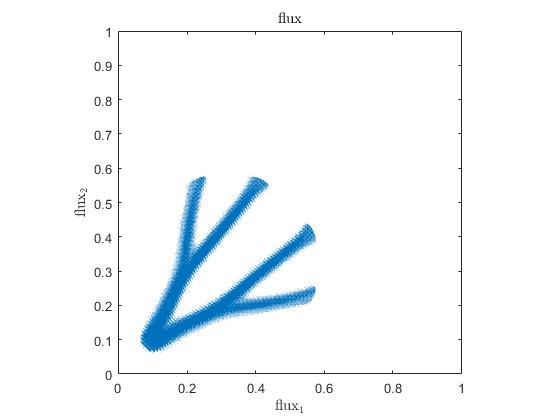}
	\end{minipage}
	\begin{minipage}{.32\textwidth}    \includegraphics[width=1\textwidth]{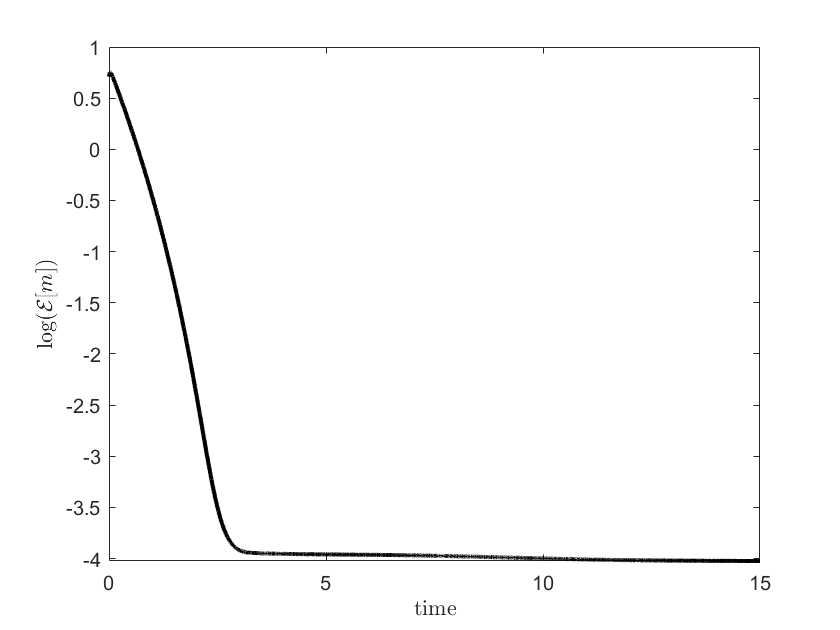}
	\end{minipage}
	\begin{minipage}{.32\textwidth}    \includegraphics[width=1.\textwidth]{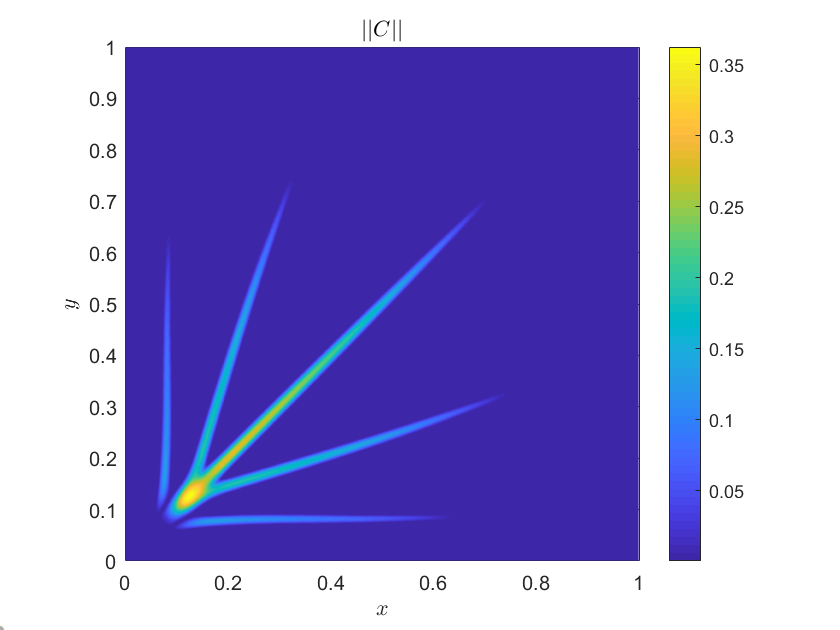}
	\end{minipage}
	\begin{minipage}{.32\textwidth}    \includegraphics[width=1\textwidth]{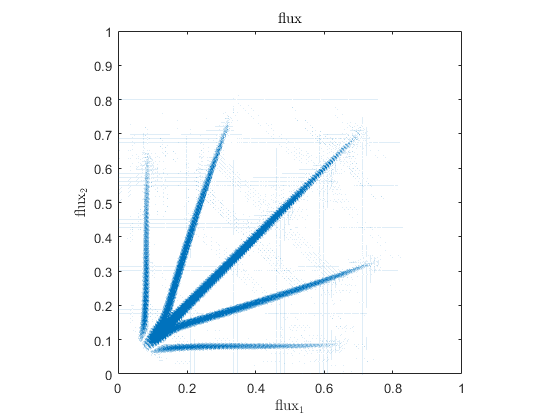}
	\end{minipage}
	\begin{minipage}{.32\textwidth}    \includegraphics[width=1\textwidth]{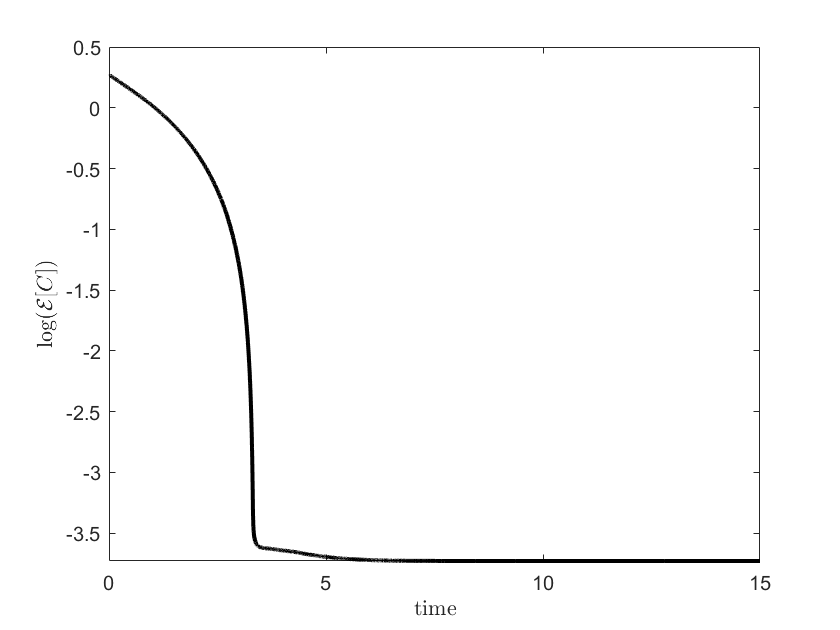}
	\end{minipage}
	\caption{\textit{In this figure we show three different quantities of the same computations, with the parameters defined in} \textsc{TestD}: \textit{the module of the variables at final time (left panels), the flux also at final time (central panels) and the energy as a function of time (right panels). The first row is about the variable $m$ and the second one is for the variable $\mathbb{C}$.}}
	\label{fig_testD}
\end{figure}

In Fig.~\ref{fig_testD} we show three different quantities of the \textsc{TestD}: the module of the variables at final time (first column), the two components of the flux
{$|\mathbb{C}\nabla p|$}  at final time (second column) and the energy as a function of time (third column). In the first row we have the results for the variable $m$ and in the second row the results for the variable $\mathbb{C}$. As expected, the energy decays in time for both variables, and it is very small at final time, which indicates that we are close to the steady state of the systems.
The main difference between the two variables are the shape of the network, with a $Y-$shape for the conductivity vector and a $V-$shape for the tensor. 

\begin{figure}[H]
	\centering
	\begin{minipage}{.32\textwidth}    \includegraphics[width=1.\textwidth]{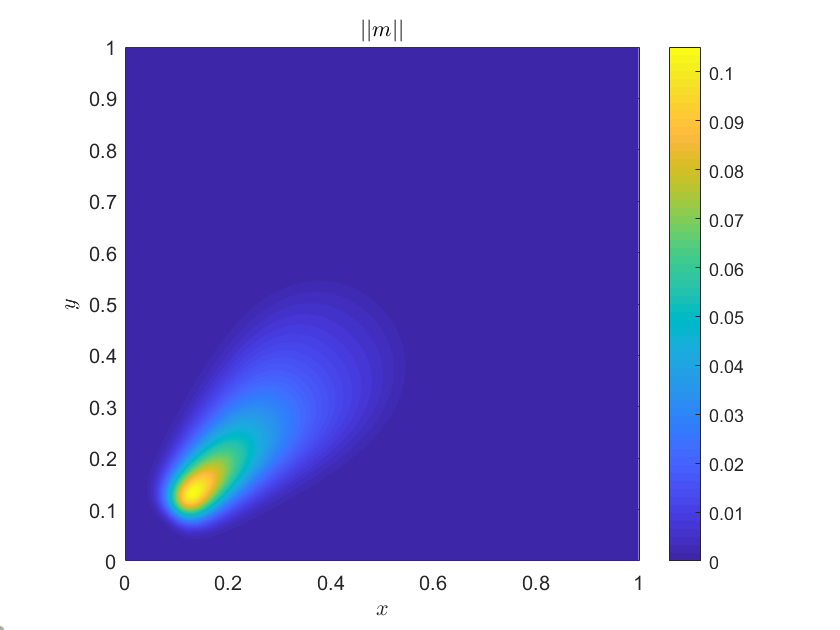}
	\end{minipage}
	\begin{minipage}{.32\textwidth}    \includegraphics[width=1\textwidth]{Figures_Africomp/TestDm_m.png}
	\end{minipage}
	\begin{minipage}{.32\textwidth}    \includegraphics[width=1\textwidth]{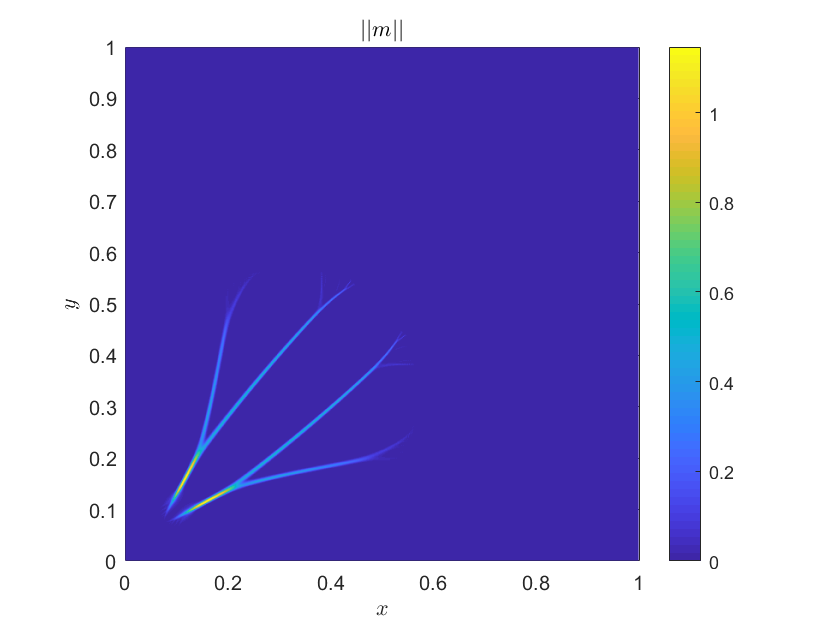}
	\end{minipage}
	\begin{minipage}{.32\textwidth}    \includegraphics[width=1.\textwidth]{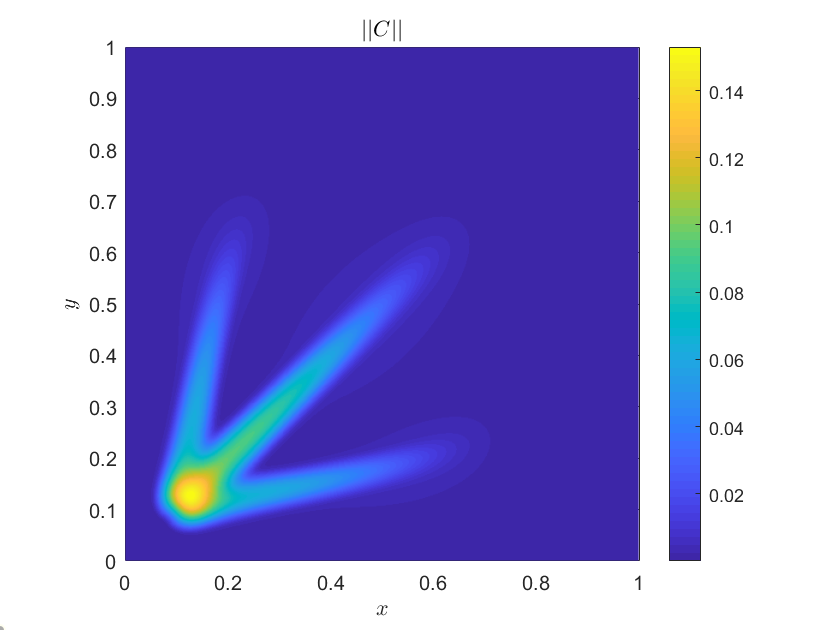}
	\end{minipage}
	\begin{minipage}{.32\textwidth}    \includegraphics[width=1\textwidth]{Figures_Africomp/C_testD.png}
	\end{minipage}
	\begin{minipage}{.32\textwidth}    \includegraphics[width=1\textwidth]{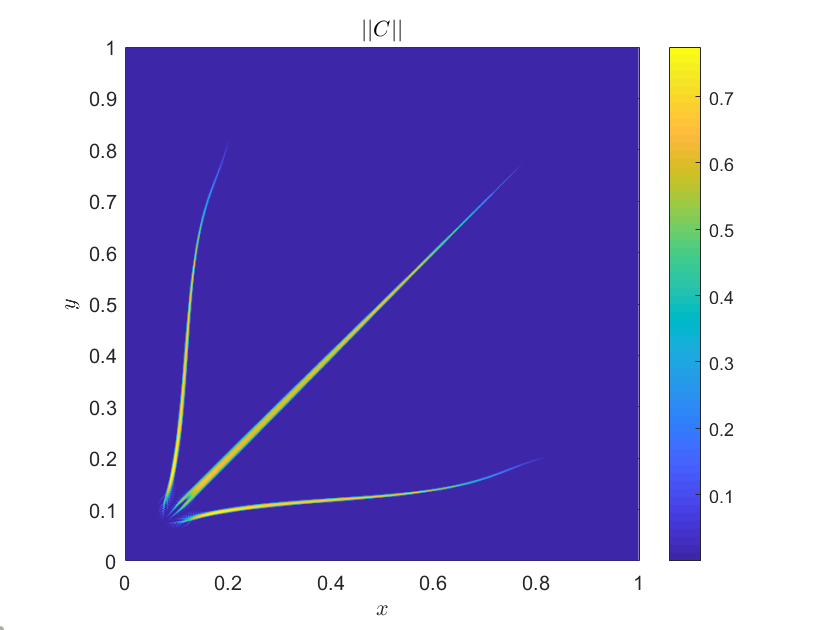}
	\end{minipage}
	\caption{\textit{In this figure we show the difference of the results on varying the diffusivity $D$. In the first column we have the results of the} \textsc{TestG} \textit{ (D = 0.05), in the second column we choose the parameters of the} \textsc{TestD} \textit{ (D = 0.01) and in the third one the } \textsc{TestE} \textit{ (D = 0.001)}. \textit{The first row shows the results of the variable $m$ and the second row those of the variable $\mathbb{C}$.}}
	\label{figure_diff}
\end{figure}

In Fig.~\ref{figure_diff} we show the results obtained when varying the parameter $D$ in Eqs.~(\ref{sy_m_p}-\ref{sy_m_diff}) and in Eqs.~(\ref{eq_darcy_p_tens}-\ref{eq_C_metab}). The tests we consider are: \textsc{TestG} (first column), \textsc{TestD} (second column) and \textsc{TestE} (third column), with $D \in \{0.05,0.01,0.001\}$, for the variable $m$ in the first row and for $\mathbb{C}$ in the second row. The ramifications become more evident when decreasing the diffusivity, and they get thinner and thinner. For the first parameter chosen, $D = 0.05$, we are not able to see those ramifications for the vector $m$ because the time-scale associated with the diffusion is too fast to capture the details.   

\begin{figure}[H]
	\centering
	\begin{minipage}{.32\textwidth}    \includegraphics[width=1.\textwidth]{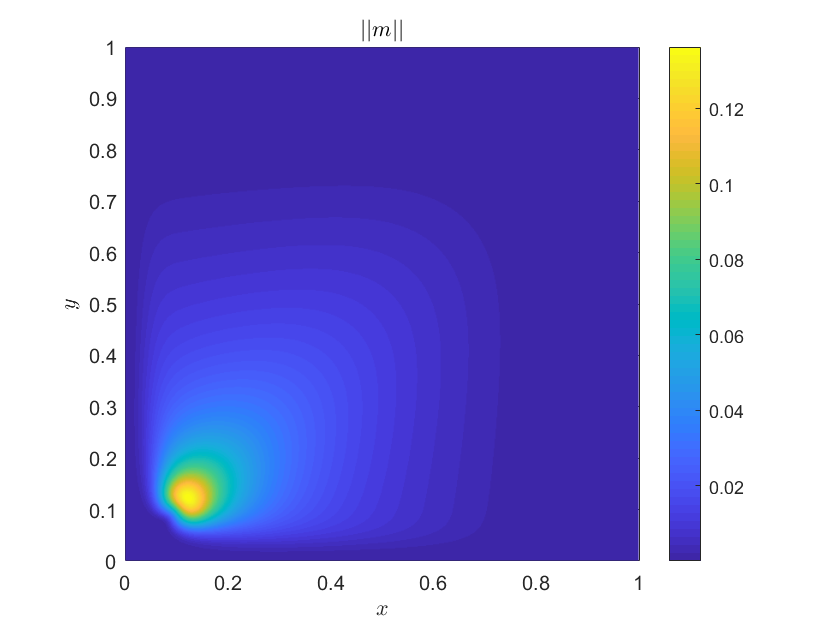}
	\end{minipage}
	\begin{minipage}{.32\textwidth}    \includegraphics[width=1\textwidth]{Figures_Africomp/TestDm_m.png}
	\end{minipage}
	\begin{minipage}{.32\textwidth}    \includegraphics[width=1\textwidth]{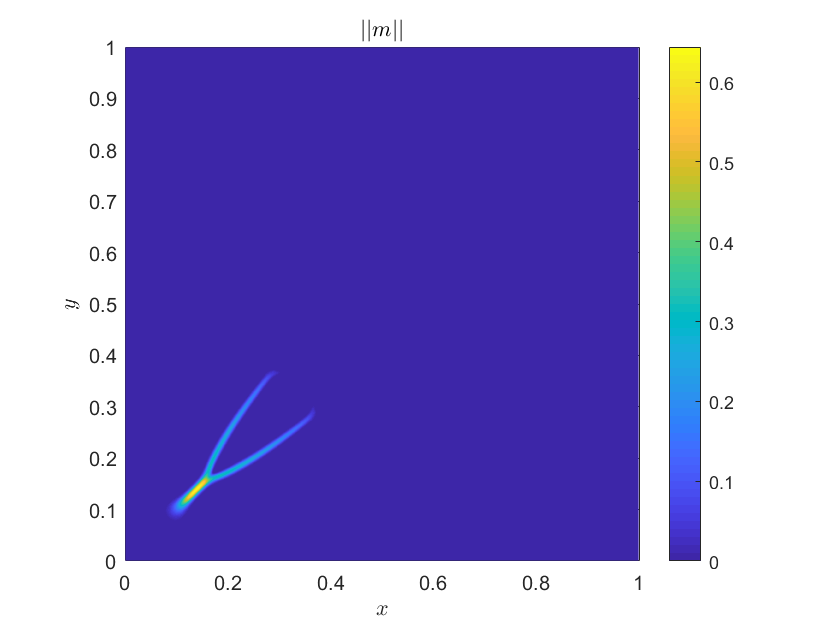}
	\end{minipage}
	\begin{minipage}{.32\textwidth}    \includegraphics[width=1.\textwidth]{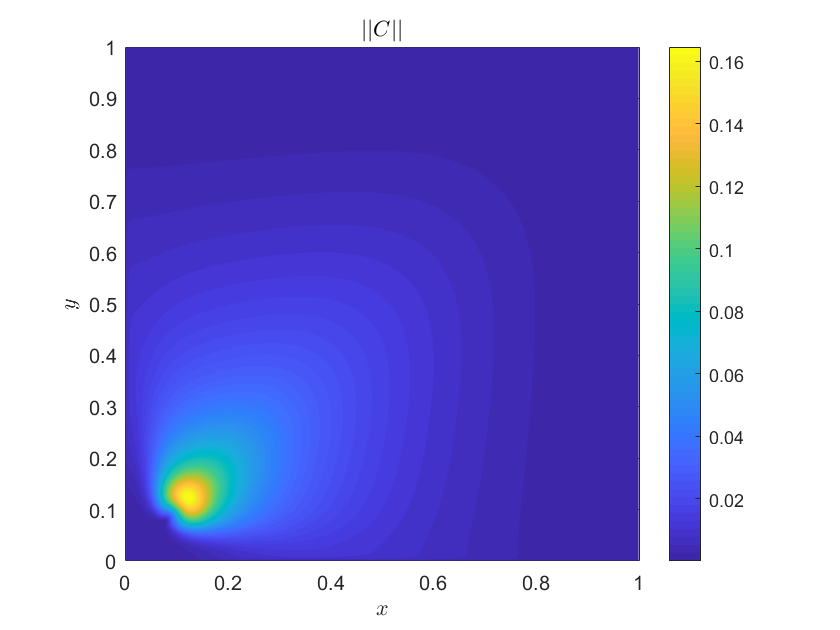}
	\end{minipage}
	\begin{minipage}{.32\textwidth}    \includegraphics[width=1\textwidth]{Figures_Africomp/C_testD.png}
	\end{minipage}
	\begin{minipage}{.32\textwidth}    \includegraphics[width=1\textwidth]{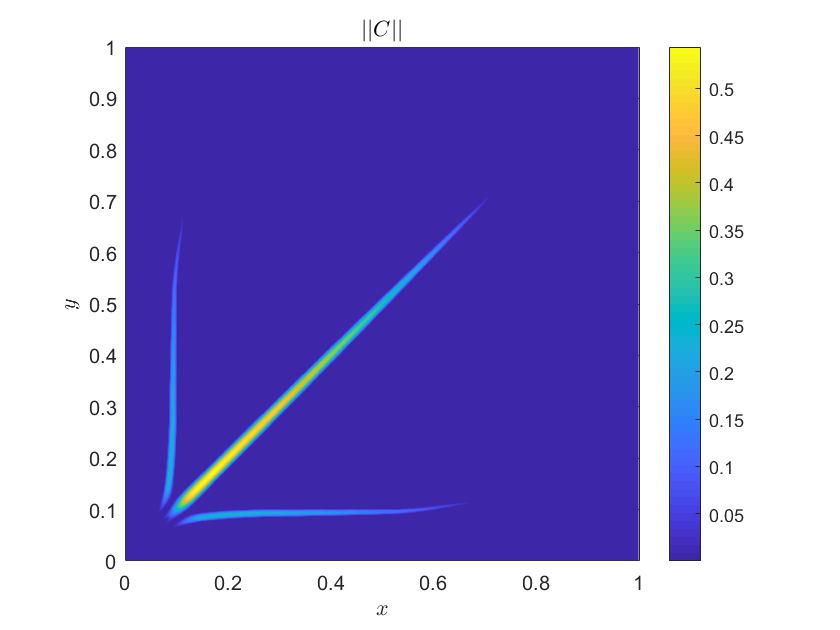}
	\end{minipage}
	\caption{\textit{In this figure we show the difference of the results on varying the relaxation exponent $\gamma$. In the first column we have the results of the} \textsc{TestH} \textit{ ($\gamma $ = 1), in the second column we choose the parameters of the} \textsc{TestD} \textit{ ($\gamma = 0.75$) and in the third one the } \textsc{TestF} \textit{ ($\gamma = 0.5$). The first row shows the results of the variable $m$ and the second row  of the variable $\mathbb{C}$.}}  
	\label{fig_gamma}
\end{figure}

In Fig.~\ref{fig_gamma} we observe the dependence on the relaxation exponent $\gamma$. In the first column we report the results of \textsc{TestH}, in the second, those corresponding to  \textsc{TestD} and in the third one, those corresponding to \textsc{TestF}. Again, in the first row we show the results for the variable $m$ and in the second row those for the variable $\mathbb{C}$. If $\gamma = 1$ the results do not show the details of the network, and it seems that $\gamma = 0.75$ is the parameter that better represents the leaf network.  

\begin{figure}[H]
	\centering
	\begin{minipage}{.32\textwidth}    \includegraphics[width=1.\textwidth]{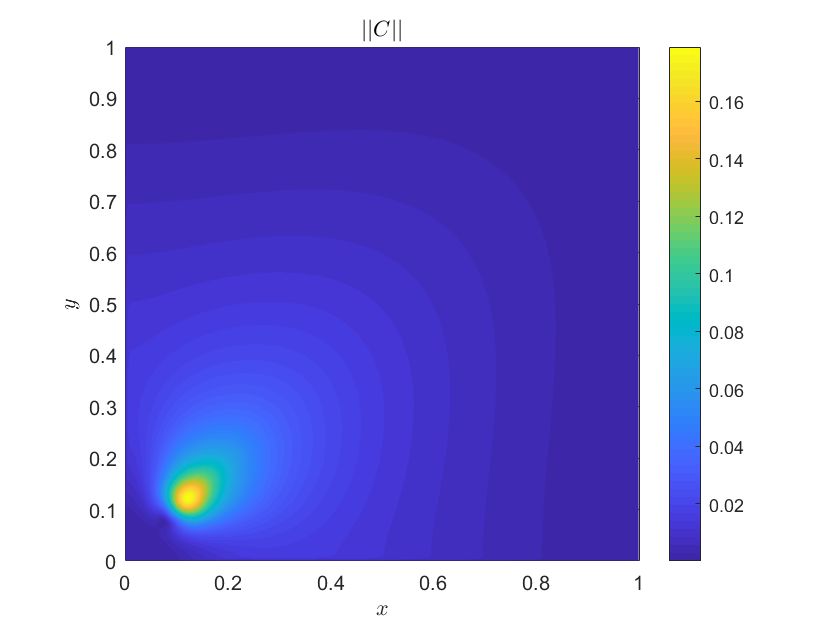}
	\end{minipage}
	\begin{minipage}{.32\textwidth}    \includegraphics[width=1\textwidth]{Figures_Africomp/C_testD.png}
	\end{minipage}
	\begin{minipage}{.32\textwidth}    \includegraphics[width=1\textwidth]{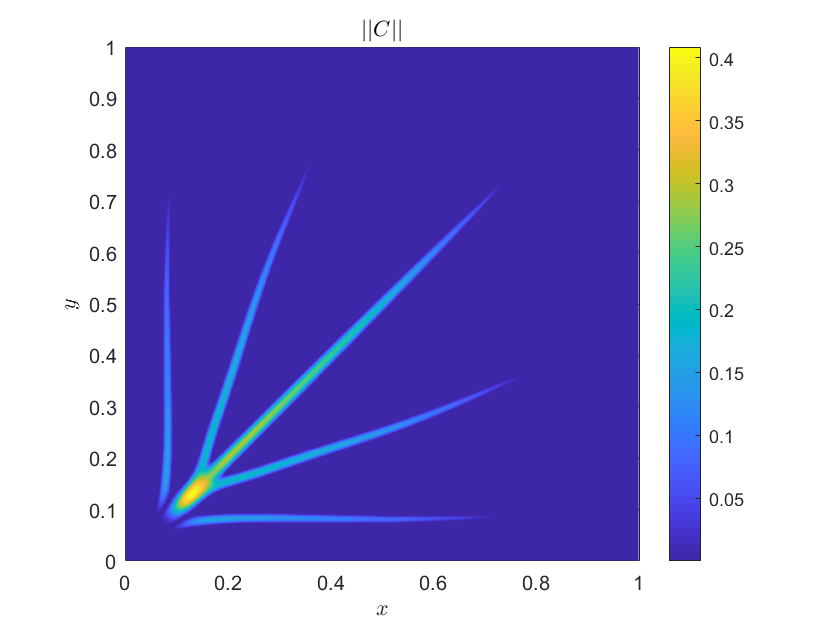}
	\end{minipage}
	\caption{\textit{In this figure we show the results for the variable $\mathbb{C}$ varying the parameter $\varepsilon$. On the left we have the results of the} \textsc{TestI} \textit{ ($\varepsilon = 10^{-2}$) in the central panel} \textsc{TestD} \textit{ $( \varepsilon = 10^{-3})$ and on the right } \textsc{TestL} $(\varepsilon = 10^{-4})$.}  
	\label{fig_epsilon}
\end{figure}

In Fig.~\ref{fig_epsilon} we show the behaviour of the solution $\mathbb{C}$, when $\varepsilon \to 0$. We see the results for $\varepsilon \in \{10^{-2} (\rm left), 10^{-3} (\rm center), 10^{-4} (\rm right)\}$, and again we notice that for the largest value of $\varepsilon$ we are not able to see any ramification. While we see that for $\varepsilon $ smaller than $10^{-3}$ we are close to the asymptotic behaviour.

\subsection{Quantitative agreement}
\label{section_quantitative}
In this section we show some quantitative comparison between the two models.
For this reason we consider \textit{well prepared} initial data and we look for compatible parameters. 

The goal of this part is to choose a convenient set of parameters in order to compare the two systems, trying to make them as close as possible. Now we distinguish the parameters $(D_l,c_l,\alpha_l,\gamma_l)$ with $l = 1$, for the $\mathbb{C}-$model, and $l=2$, for the $m-$model.

For simplicity, the choice of parameter is performed by comparing the two models in  one space dimension.
In 1D the systems~(\ref{eq_darcy_p_tens},\ref{eq_C_metab}) and (\ref{sy_m_p}-\ref{sy_m_diff}) read
\begin{align}
	\label{eq_C_1D}
	C_t - D_1^2\,C_{xx} - c_1^2\,p_x^2 + \alpha_1|C|^{\gamma_1 - 2} C &= 0 \\
	\label{eq_m_1D}
	m_t - D_2^2\,m_{xx} - c_2^2\,p_x^2\, m + \alpha_2|m|^{2(\gamma_2 - 1)} m&= 0. 
\end{align}

Now we suppose that $C$ has the following form
\begin{equation}
	\label{eq_expr_C_B}
	C = m^2 + B,
\end{equation}
where $B$ is a measure of the discrepancy between the two models, and we set the initial conditions so that $B(t=0) = 0$. If we substitute  Eq.~\eqref{eq_expr_C_B} in Eq.~\eqref{eq_C_1D}, we have
\begin{align}
	\label{eq_m_sq_1}
	(m^2)_t - D_1^2(m^2)_{xx} - c_1^2p_x^2 + \alpha_1|m^2 + B|^{\gamma - 2} m^2 = -B_t + D_1^2 B_{xx} - \alpha_1 |m^2 + B|^{\gamma - 2} B.
\end{align}
At this point we multiply Eq.~\eqref{eq_m_1D} by a factor $(2 m)$, and we obtain
\begin{equation}
	\label{eq_m_sq_2}
	2m\, m_t - 2 D_2^2 m\, m_{xx} - 2\, c_2^2\, p_x^2\, m^2 + 2\alpha_2 |m|^{2(\gamma_2-1)}m^2 = 0.
\end{equation}
After some manipulation, Eqs.~(\ref{eq_m_sq_1},\ref{eq_m_sq_2}) become:
\begin{eqnarray}
	\label{eq_manip_C}
	&2 m\, m_t - D_1^2 (m^2)_{xx} - c_1^2p_x^2 + \alpha_1\,(m^2)^{\gamma - 1} = \underbrace{-B_t + D_1^2 B_{xx} - \alpha_1 |m^2|^{\gamma - 2} B}_{:= \mathcal{R}} &\\ \label{eq_manip_m}
	&2m\, m_t - D_2^2 (m^2)_{xx} + 2D_2^2 (m_x)^2 - 2\, c_2^2\, p_x^2\, m^2 + 2\alpha_2 (m^2)^{\gamma_2} = 0 &
\end{eqnarray}
where $\mathcal{R}$ is the residual. We made use of the following identity in the equation for $m$
\begin{equation}
	\label{eq_identity}
	(m)^2_{xx} = 2 m\,m_{xx} + 2\left( (m_x)^2 \right).
\end{equation} 
Now we consider the difference between Eq.~\eqref{eq_manip_C} and Eq.~\eqref{eq_manip_m}, and we obtain
\begin{equation}
	(D_2^2 - D_1^2) (m^2)_{xx} -  2D_2^2 (m_x)^2 + ( 2\, c_2^2 m^2 - c_1^2) p_x^2 + \alpha_1\,(m^2)^{\gamma - 1} - 2\alpha_2 (m^2)^{\gamma_2} = \mathcal{R}
\end{equation}
Since we want the residual to be small, in absolute value, a convenient choice for the sets of variables is the following
\begin{equation}
	\label{eq_set_of_param}
	D_1 = D_2, \quad \alpha_1 = 2\alpha_2, \quad \gamma_1 = \gamma_2 + 1, \quad c_1 = \sqrt{2} \, c_2 |m|
\end{equation} 
and for the initial conditions we choose $m^0$ and $C^0$ such that, initially, we have
\[ 
C^0 = (m^0)^2 = {\rm const}.
\]
In this way the first derivative in space $m_x$ is also equal to 0 after one time step.

In order to show some results in 2D, we need to define the initial conditions for $m^0_{\rm comp}$ and $\mathbb{C}^0_{\rm comp}$ such that, at initial time, we have again 
\[ 
\mathbb{C}^0 = m^0 \otimes m^0,
\]
with $m^0_{\rm comp} = [\sqrt{2}/2,\sqrt{2}/2]^T$ and $\mathbb{C}^0_{\rm comp } = [0.5,0.5,0.5]^T$, while the values of the parameters are defined in \textsc{TestM}. In 2D, the equivalent expression of \eqref{eq_expr_C_B} is
\[ \mathbb{C} = m\otimes m + \mathbb{B}.
\]

In this subsection we comment on the solutions of the following tests

\begin{table}[H]
	\centering
	\begin{tabular}{|c|c|c|c|c|c|c|c|c|}  \hline
		& & $\alpha_1, \alpha_2$ & $c_1, c_2$ & $D_1 = D_2$ & $\varepsilon$ & $\gamma_1, \gamma_2$ & $r$ & $N$ \\ \hline
		set of parameters & \textsc{TestM}:&  1, 0.5 & $\sqrt{2}$, 1 & 0.1 &  $ 10^{-1}$ & 1.75, 0.75 & 0.1 & 600 \\ \hline
	\end{tabular}
	\caption{\textit{In this table we define the two set of parameters in Eq.~\eqref{eq_set_of_param}. }}
	\label{tab_set_param}
\end{table}        

\begin{table}[H]
	\centering
	\begin{tabular}{|c|c|c|c|c|c|c|c|c|c|}  \hline
		& & $\alpha$ & $c$ & $D$ & $\varepsilon$ & $\gamma$ & $r$ & $t_{\rm fin}$ & $N$ \\  \hline
		$D = 0$ & \textsc{TestN}: &  0.75 & 5 & 0 & $10^{-3}$ & 0.75 & 0.005 & 15 & 600 \\ \hline
		$D = 10^{-5}$ & \textsc{TestO}: &  0.75 & 5 & $10^{-5}$ & $10^{-3}$ & 0.75 & 0.005 & 15 & 600 \\ \hline
	\end{tabular}
	\caption{\textit{In this table we define the two tests showed in Figs.\ref{fig_zero_diffusivity}-\ref{fig_zero_diffusivity_2}, with $D = 0 $ (\textsc{TestN}) and $D = 10^{-5}$ (\textsc{TestO}). }}
	\label{tab_zero_diff}
\end{table} 

In Table~\ref{table_two_sol} we show the time evolution of the norm of the difference between $\mathbb{C}$ and $m\otimes m$, to see how they move away from each other when we increase the time. In this table we see that the two solutions are very different, even after few time steps, and for this reason, they are difficult to be compared. The definition of  $||\mathbb{B}||$ is the following
\[ ||\mathbb{B}|| := \frac{\Big|\Big| |\mathbb{C}| - | m \otimes m| \Big|\Big|}{\Big|\Big| |m \otimes m| \Big|\Big|}
\]
after $n_t$ time steps, with ${\rm time} = n_t\Delta t$, such that $n_t = 6(2^k),\, k = 0,1,2,3,4,5,6$, 
where $|\mathbb{C}| := \sqrt{\mathbb{C}_{11}^2 + 2\,\mathbb{C}_{12}^2 + \mathbb{C}_{22}^2}$ and \\$|m \otimes m| := \sqrt{(m_{1}^2)^2 + 2\,(m_1\,m_2)^2 + (m_{2}^2)^2}$. 
\begin{table}[H]
	\centering
	\begin{tabular}{||c||c||c||c||c||c||c||c||}
		\hline \hline
		time & 0.01 & 0.02 & 0.04 & 0.08 & 0.16 & 0.32 & 0.64  \\
		\hline \hline
		$||\mathbb{B}||$ & 0.0348  &  0.0538  &  0.0697  &  0.0982  &  0.1509 &  0.2611   & 0.5320 \\
		\hline \hline     
	\end{tabular}
	\caption{\textit{Here we show the values of $\mathbb{B}$, after $n_t$ time steps, where $n_t = 2^k, k = 0,1,2,3,4,5,6$.}}
	\label{table_two_sol}
\end{table}

{As it appears from the table, the two models move quickly far apart from each other, suggesting an intrinsically different behaviour.}

\begin{figure}[H]
	\centering
	\begin{minipage}{.49\textwidth}    \includegraphics[width=1.\textwidth]{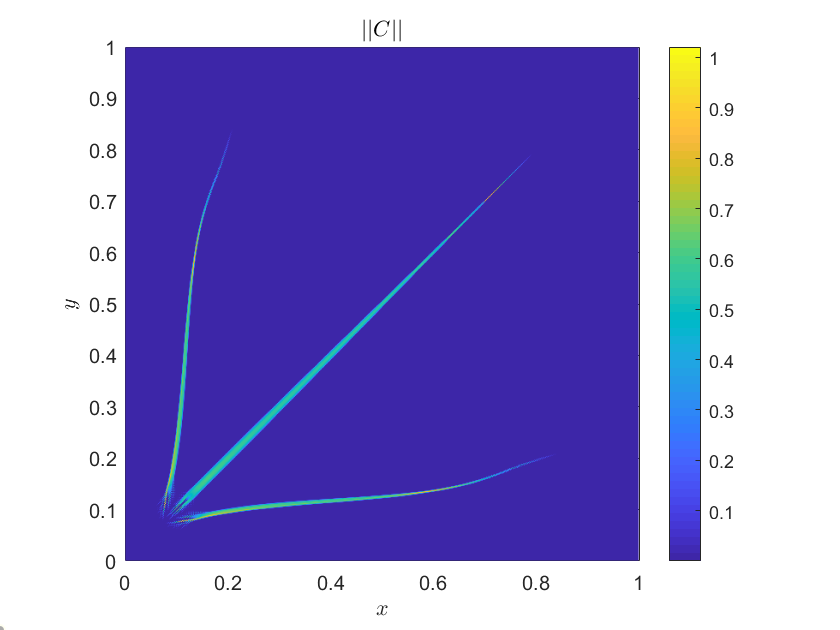}
	\end{minipage}
	\begin{minipage}{.49\textwidth}    \includegraphics[width=1\textwidth]{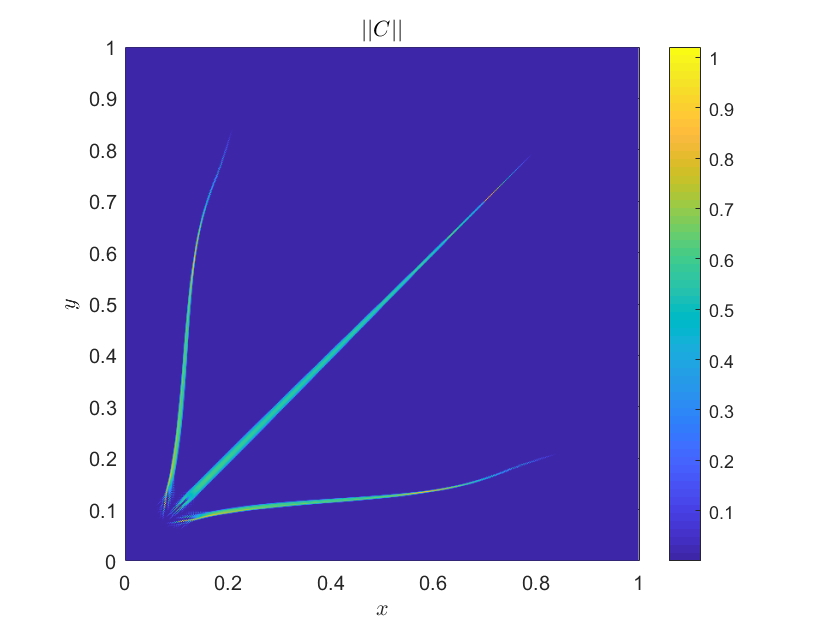}
	\end{minipage}
	\caption{\textit{Comparison between \textsc{TestN} and \textsc{TestO}, with initial condition defined in Eq.~\eqref{eq_ic_m}. The main difference is that, on the left we have $D = 0$, while, on the right, $D = 10^{-5}$.}}  
	\label{fig_zero_diffusivity}
\end{figure}

Now we are interested in showing the solution of the $\mathbb{C}-$model, in the case of zero-diffusivity. Since the randomness of the network is common in nature but is also very effective in stabilizing the equations, we want to see if there is some analogy in considering the cases $D = 0$ and $D \ll 1$. In Fig.~\ref{fig_zero_diffusivity} we show the agreement of the two solutions, with $D = 0$ (left panel) and $D = 10^{-5}$ (right panel). The other parameters are defined in \textsc{TestN} and \textsc{TestO}, while the initial condition is defined in Eq.~\eqref{eq_ic_m} for Fig.~\ref{fig_zero_diffusivity}.

\begin{figure}[H]
	\centering
	\begin{minipage}{.49\textwidth}    \includegraphics[width=1.\textwidth]{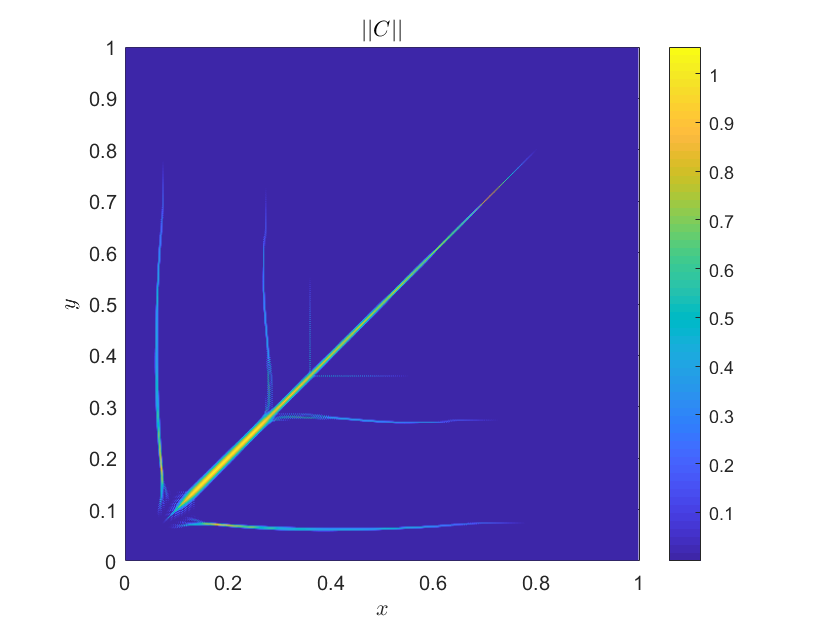}
	\end{minipage}
	\begin{minipage}{.49\textwidth}    \includegraphics[width=1\textwidth]{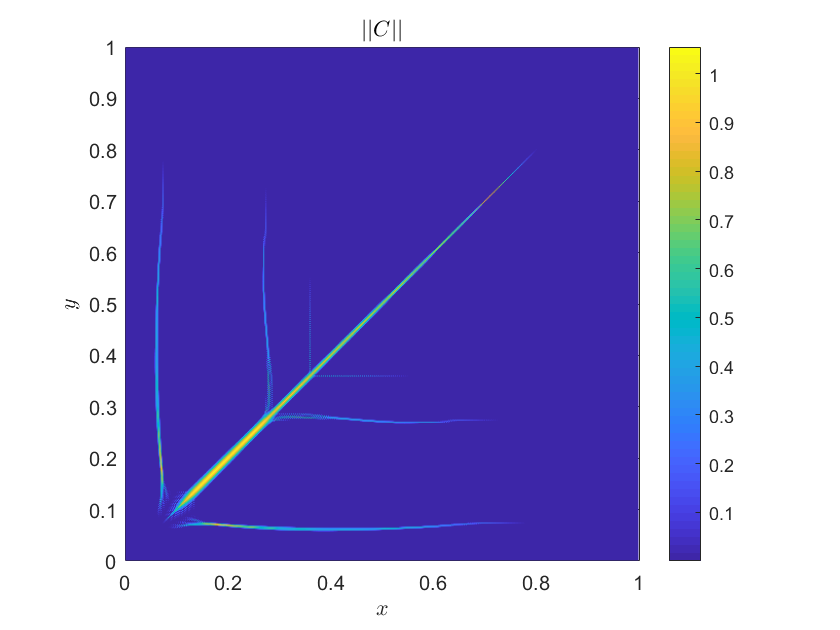}
	\end{minipage}
	\caption{\textit{Comparison between \textsc{TestN} and \textsc{TestO}. Here the initial condition is defined in Eq.~\eqref{eq_ic_2}, with $D = 0$ on the left and $D = 10^{-5}$ on the right.}}  
	\label{fig_zero_diffusivity_2}
\end{figure}

In Fig.~\ref{fig_zero_diffusivity_2} we {illustrate the long time solution for the $\mathbb{C}$-model obtained with the following space dependent initial condition}:
\begin{equation}
	\label{eq_ic_2}
	\mathbb{C}^0 = f(x,y) \mathbb{I} , \quad f(x,y) = (2 - |X + Y|)\exp(-10(|X-Y|))
\end{equation}

We also note that, if we set the diffusion coefficient equal to zero, the model for $m$ reduces to a reaction equation for the conductivity. This means that {if $\gamma>1/2$ the support of the unknown $m$ remains unchanged. In particular, it cannot extend, while in some regions the numerical support (i.e.\ the region in which $|m|$ becomes lower than a given small threshold) may shrink.} 

\paragraph{Alternative boundary conditions}
{Furthermore, in the case of zero-diffusivity for the $m-$model, we observe an anomalous behavior of the solution near the boundaries. In order to overcome such a problem we propose  an \textit{ad hoc} boundary condition as illustrated below.}

Let us consider the equations for $m$ with $D = 0$ and in the limit of steady state. We have
\begin{align*}
	c^2\nabla p\otimes \nabla p \,m -\alpha |m|^{2(\gamma - 1)} m & = 0 \quad {\rm in} \, \partial \Omega  
\end{align*}
that we can write as follows
\begin{align}
	\label{eq_m_gradp}
	(m \cdot \nabla p)\nabla p = \frac{\alpha}{c^2} |m|^{2(\gamma - 1)} m \quad {\rm in} \, \partial \Omega  
\end{align}
Thus, we can deduce that $m\propto\nabla p$.  This means that there exists a constant $\beta$, such that, 
\begin{equation}
	\label{eq_m_prop_p}
	m = \beta \nabla p.
\end{equation}
If we substitute the Eq.~\eqref{eq_m_prop_p} in Eq.~\eqref{eq_m_gradp}, we obtain
\begin{equation}
	\beta |\nabla p|^2 \nabla p = \frac{\alpha}{c^2} \beta^{2(\gamma - 1)}|\nabla p|^{2(\gamma - 1)} \beta \nabla p
\end{equation}
and if we solve it for $\beta$, we have the boundary condition for $m$
\begin{align}
	\label{eq_m_p_bc}
	\left. m \right|_{\partial \Omega} = \underbrace{\left(\frac{c^2}{\alpha} |\nabla p|^{4 - 2\gamma}  \right)^{\frac{1}{2(\gamma - 1)}}}_{\beta} \nabla p.
\end{align}

Analogously, we can find also a boundary condition for the vector $\mathbb{C}$. Again, starting with zero-diffusivity and at steady state, we have
\begin{equation}
	\label{eq_C_gradp}
	\nabla p \otimes \nabla p = \frac{\alpha}{c^2} |\mathbb{C}|^{\gamma - 2}\mathbb{C} \quad {\rm in} \, \partial \Omega
\end{equation}
that means that the conductivity is proportional to the tensor product of the pressure gradient, i.e. $\mathbb{C} \propto \nabla p \otimes \nabla p$. Again, we look for a constant $\tilde{\beta}$, such that,
\begin{equation}
	\label{eq_C_p}
	\mathbb{C} = \tilde{\beta} \nabla p \otimes \nabla p.
\end{equation}
Now we substitute the Eq.~\eqref{eq_C_p} in Eq.~\eqref{eq_C_gradp}, and we solve for $\tilde{\beta}$. In this way, as before, we find the expression for the conductivity tensor at the boundary
\begin{equation}
	\label{eq_C_p_bc}
	\left. \mathbb{C}\right|_{\partial \Omega} = \underbrace{\left( \frac{c^2}{\alpha} |\nabla p|^{-2(\gamma - 2)} \right) ^{\frac{1}{\gamma - 1}}}_{\tilde{\beta}}\nabla p \otimes \nabla p.
\end{equation}

Conditions (\ref{eq_m_p_bc},\ref{eq_C_p_bc}) might be a reasonable choice in the case of zero diffusivity. This treatment has the drawback of introducing additional non-linearity to the system.
Alternative boundary conditions are currently under investigation.

\begin{figure}[H]
	\centering
	\begin{minipage}{.49\textwidth}    \includegraphics[width=1.\textwidth]{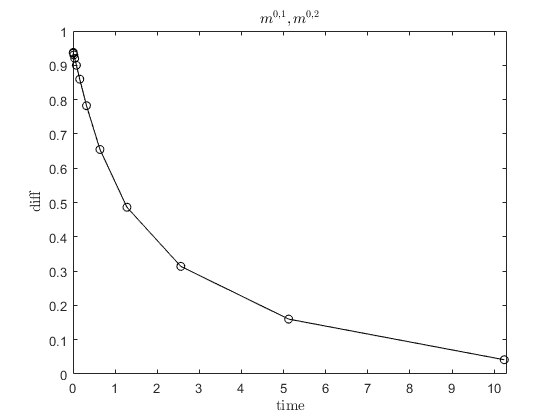}
	\end{minipage}
	\begin{minipage}{.49\textwidth}    \includegraphics[width=1\textwidth]{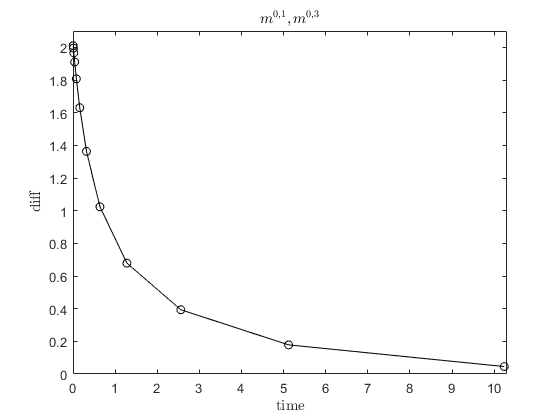}
	\end{minipage}
	\caption{\textit{In this figure we show the difference between two different solutions choosing, as initial condition, $m^{0,1}, m^{0,2}$ (on the left) and $m^{0,1}, m^{0,3}$ (on the right). We plot the expression defined in Eq.~\eqref{eq_diff}, as a function of time, with $\gamma = 1.75 > 1$, and the others parameters are defined in \textsc{TestD}. }}
	\label{fig_steady_state_gg}
\end{figure}

Another aspect we want to focus on, is the \textit{steady state} for the $m-$model. In two different cases: $\gamma < 1$ and $\gamma > 1$ (as the authors show in \cite{marko_perthame_2}). For this reason we define different initial conditions for the vector $m$, such that
\begin{eqnarray*}
	&  m^{0,1}_1 = 1, \quad m^{0,1}_2 = \sqrt{2};\quad 
	m^{0,2}_1 = 5, \quad m^{0,2}_2 = 5; & \\
	& m^{0,3}_1 = (2 - |X + Y|)\exp{(-10|X-Y|)}, \quad m^{0,3}_2 = m^{0,3}_1;& 
\end{eqnarray*} 
and in Fig.~\ref{fig_steady_state_gg} we plot the following quantity
\begin{equation}
	\label{eq_diff}
	{\rm diff} = \frac{|| m^\beta - m^\eta ||}{|| m^\eta ||}
\end{equation}
where $||\cdot ||$ is the Frobenius norm, and $\beta,\eta \in \{1,2,3\}$. In this way we see the difference between the solutions (with initial condition $m^{0,1}$ and $m^{0,2}$ in the left panel and with $m^{0,1}$ and $m^{0,3}$ in the right panel), as function of time. In this way, we support with numerical evidences that, for the $m-$model and for $\gamma > 1$, the steady state is unique and it does not depend on the initial conditions, as expected (see e.g. \cite{marko_perthame_2}). 

\begin{figure}[H]
	\centering
	\begin{minipage}{.49\textwidth}    \includegraphics[width=1\textwidth]{Figures_Africomp/TestDm_m.png}
	\end{minipage}
	\begin{minipage}{.49\textwidth}    \includegraphics[width=1\textwidth]{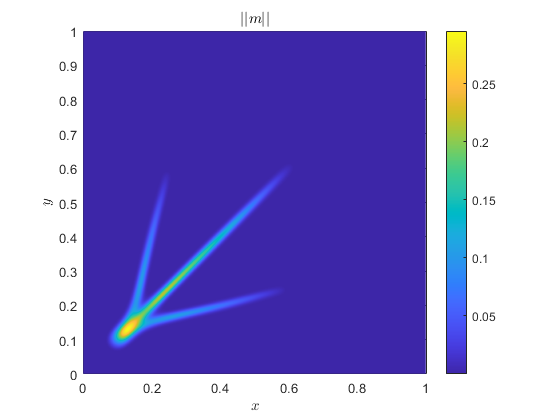}
	\end{minipage}
	\caption{\textit{In this figure we show the steady states when $\gamma = 0.75 < 1$. On the left the initial condition is $m^0 = m^{0,1} = 1$ while, on the right, the initial condition is a function of space, $m^{0} = m^{0,3}$, and the parameters are defined in \textsc{TestD}.}}
	\label{fig_steady_state_gl}
\end{figure}

In Fig.~\ref{fig_steady_state_gl}, that is the case for $\gamma < 1$, we see that we reach two different steady states, when choosing two different initial conditions, {suggesting that the steady state solution is not unique when $\gamma<1$.}

\section{Conclusions}
In this paper we use an elliptic-parabolic model to study the formation of biological network, and, in particular, of leaf venation networks. Throughout the paper we compare the solutions of two different systems, that derive from the Cai-Hu model, one for the conductivity vector and one for the conductivity tensor, and we explore the dependence of the solution on the parameters of the two models. 

As we said before, all the components of the two systems are very stiff. In particularly, the $\mathbb{C}$-system is more challenging because of the negative exponent in the reaction term. For this reason, we add a \textit{regularization parameter} $\varepsilon$, and we compute numerical solutions for smaller and smaller values of $\varepsilon$. This parameter prevents the instability coming from the division by zero in the reaction term.

We make use of finite differences scheme to compute the two solutions, with central differences for the space discretization and a \textit{symmetric-}ADI method in time. 
The convergence rate is calculated numerically and denotes the second order accuracy. 

At the end we added some quantitative comparison between the two systems, choosing more suitable sets of parameters. We see that the two solutions differ significantly, even after few time steps. This aspect makes it problematic any kind of direct comparison between the two systems.

Then, we showed some results in the case of zero diffusivity for the conductivity tensor. The last tests we consider are in agreement with the results achieved in \cite{marko_perthame_2}, where the authors prove that there is a unique steady state for the $m-$system when the metabolic exponent $\gamma$ is greater than 1, {and provide evidence that this is not true when $\gamma<1$.}

\section*{Acknowledgments}
G.R. thanks ITN-ETN Horizon 2020 Project ModCompShock, Modeling and Computation on Shocks and Interfaces, Project Reference 642768, and the Italian Ministry of Instruction, University and Research (MIUR) to support this research with funds coming from PRIN Project 2017 (No.2017KKJP4X entitled ''Innovative numerical methods for evolutionary partial differential equations and applications''.

\bibliographystyle{unsrt}
\bibliography{main}

\end{document}